\documentclass[12pt,a4paper]{amsart}
\usepackage[utf8]{inputenc}
\usepackage{tikz}
\usepackage[english]{babel}
\usepackage[T1]{fontenc}
\usepackage{amsmath}
\usepackage{amsfonts}
\usepackage{amssymb}
\usepackage{graphicx}
\usepackage[left=2.5cm,right=2.5cm,top=3cm,bottom=2.5cm]{geometry}
\usepackage{color}

\theoremstyle{definition}
\newtheorem{defi}{Definition}[section]
\newtheorem{ex}[defi]{Example}
\newtheorem{rem}[defi]{Remark}
\newtheorem{rems}[defi]{Remarks}

\theoremstyle{plain}
\newtheorem{thm}[defi]{Theorem}
\newtheorem{prop}[defi]{Proposition}
\newtheorem{lem}[defi]{Lemma}
\newtheorem{cor}[defi]{Corollary}
\newtheorem{conj}[defi]{Conjecture}

\numberwithin{equation}{section}

\def\R{{\mathbb{R}}}
\def\N{{\mathbb{N}}}
\DeclareMathOperator{\dist}{dist}
\newcommand{\dx}{{\rm d}x}
\newcommand{\dt}{{\rm d}t}

\begin{document}

\title[The fundamental gap of quantum graphs]{Optimizing the Fundamental Eigenvalue Gap of Quantum Graphs}

\author[M. Ahrami]{Mohammed Ahrami}
\author[Z. El Allali]{Zakaria El Allali}
\author[E.M. Harrell II]{Evans M. Harrell II}
\author[J.B. Kennedy]{James B. Kennedy}

\address{Modeling and Scientific Computing Team, 
Multidisciplinary Faculty of Nador, Mohammed I University, Morocco}
\email{m.ahrami@ump.ac.ma}

\address{Modeling and Scientific Computing Team, 
Multidisciplinary Faculty of Nador, Mohammed I University, Morocco}
\email{z.elallali@ump.ma}

\address{School of Mathematics,
Georgia Institute of Technology,\
Atlanta GA 30332-0160, USA} 
\email{harrell@math.gatech.edu}

\address{Grupo de F\'isica Matem\'atica {\rm and} Departamento de Matem\'atica, Faculdade de Ci\^encias da Universidade de Lisboa, 1749-016 Campo Grande, Lisboa, Portugal}
\email{jbkennedy@ciencias.ulisboa.pt}



\thanks{Some parts of this work (in particular many of the existence and characterization results in Sections~\ref{sec:existence} and~\ref{sec:characterization}, respectively) are evolutions of ideas in an unpublished manuscript of the third author together with Joachim Stubbe and Timo Weidl, whose role as intellectual godparents of the results is explicitly acknowledged. The work of M.A. and J.B.K. was supported by the Funda\c{c}\~ao para a Ci\^encia e a Tecnologia, Portugal, within the scope of the project NoDES: new horizons in regularity, dynamics and asymptotic analysis, reference PTDC/MAT-PUR/1788/2020, and grant UIDB/00208/2020.}

\maketitle

\begin{abstract}
We study the problem of minimizing or maximizing the fundamental spectral gap of Schr\"odinger operators on metric graphs with
either a convex potential or 
a ``single-well'' potential on an appropriate specified subset. (In the case of metric trees, such a subset can be the entire graph.)  In the convex case we find
that the minimizing and maximizing potentials are piecewise linear with only a finite number of points of non-smoothness, but give examples showing that the optimal potentials need not be constant. This is a significant departure from the usual scenarios on intervals and domains where the constant potential is typically minimizing. In the single-well case we show that the optimal potentials are piecewise constant with a finite number of jumps, and in both cases give an explicit estimate on the number of points of non-smoothness, respectively jumps, the minimizing potential can have. Furthermore, we show that, unlike on domains, it is not generally possible to find nontrivial bounds on the fundamental gap in terms of the diameter of the graph alone, within the given classes.
\end{abstract}

\section{Introduction}

This paper delves into the problem of minimizing or maximizing the fundamental spectral gap of Schr\"odinger operators on metric graphs $G$
under certain constraints on the potential energy.  Specifically, we consider potentials that are either convex or of single-well form on an appropriately specified subset. 

The fundamental spectral gap, a pivotal concept in quantum mechanics which represents the energy difference between the ground state and the first excited state of a quantum system, is defined as
$$ \Gamma[q]:=\lambda_{2}(q)- \lambda_{1}(q), $$
where $q$ is the potential and $\lambda_1(q) < \lambda_2(q)$ are the first two eigenvalues of the associated Schr\"odinger operator.

There has been extensive work on this problem in the past. In one dimension, on an interval of length $D$, the long-standing conjecture that the constant potential minimizes the fundamental gap among
all convex potentials and satisfies the lower bound  
\begin{equation}\label{conjecture}
 \Gamma[q]\geqslant \frac{3\pi^{2}}{D^{2}}, 
\end{equation}
was proved by Lavine \cite{Lavine} in $1994$. In outline, Lavine used methods of ordinary differential equations
to study the normalized eigenfunctions $u_1 \sim \lambda_1 (q)$, $u_2 \sim \lambda_2 (q)$, in particular to understand the derivative of the ratio $\frac{u_{2}}{u_{1}}$ and the subintervals where $u_{2}^{2}-u_{1}^{2}$ is positive or negative, drawing on ideas of Ashbaugh and Benguria \cite{Ash1}. He then applied a perturbative argument to show that the case of minimality is attained within the subclass of convex $q(x)$ of the form $ax+b$, and finally analyzed this special case. Some time later, in \cite{Andrews}, Andrews and Clutterbuck showed that the
fundamental gap of a Dirichlet Schr{\"o}dinger operator with semi-convex potential $q$ on a bounded convex domain $\Omega\in \mathbb{R}^{d}$ with diameter $D$ satisfies the estimate \eqref{conjecture}, thus proving the so-called fundamental gap conjecture in essentially full generality. 

Our goal is to understand the extent to which such lower bounds on the spectral gap hold, or indeed make sense, in the context of quantum graphs, that is, Schr{\"o}dinger-type operators defined on compact metric graphs $G$.

Estimates on the spectral gap of quantum-graph Laplacians in function of geometric properties of the graph (including its diameter) have been heavily studied in the last ten years or so; see, e.g., \cite{Berkolaiko-1,borthwick,KurasovBook} and the references therein. 

This paper explores the interpretation of the concepts of convexity and single-well functions within the context of a metric graph. 
We bring specific techniques to analyze the problem of optimizing the fundamental spectral gap of quantum graphs, providing insights into how different classes of potentials impact the spectral properties of quantum graphs.
The structure of the graph becomes significant when we delve into the complexities of numerous intersecting paths and the consistency conditions imposed on the function at these intersections. It is straightforward to define global convexity if the graph is a tree, but it seems rather limiting when the
graph contains a cycle: If convexity is imposed on any subgraph containing the cycle, then the function is forced to be constant on the cycle, taking on its minimum value
there. Imposing convexity on only some paths leaves us with a richer family of spectral problems, to which our methods apply. In fact always one can always define a family of functions that are globally convex on a metric graph in the sense that then functions are globally convex on a tree that includes all of the edges of the graph. Granted, there is some arbitrariness in the choice of such maximal trees. Interestingly, the original definition of single-well function on an interval does not extend on metric graph  even for trees, because different paths might for example be disjoint. This highlights the unique challenges and considerations when extending concepts from intervals to more complex structures like metric graphs. 

The organization of this paper is as follows. 
Notation and background results related to the fundamental spectral
gap of Schr\"odinger operators on quantum graphs are presented in \S \ref{sec:Prelims}, including a precise definition of the two main classes of potentials we will consider, convex and ``single-well,'' in the context of metric graphs (\S \ref{sec:convex-single-well}). There we also (\S \ref{sec:feynman-hellman}) derive a perturbation formula of Kato (or Feynman--Hellman) type (Proposition~\ref{degpert}/Corollary~\ref{cor:feynman-hellman}) for the eigenvalues, which gives general necessary conditions for the potential energy to be an optimizer (either maximizer of minimizer), and which will later be used to characterize them through arguments by contradiction. This is explicitly and necessarily formulated to be valid in the case of a degenerate eigenvalue; while folklore, such results are not easy to find in the literature. In \S \ref{sec:existence} we then establish (Theorem \ref{thm:exist}) general existence results by invoking compactness of the sets of functions with the properties of interest, which in turn follow from a general argument using either the Blaschke selection principle or the Helly compactness theorem.

With these results in hand we turn to characterizing the
optimizers in \S \ref{sec:characterization}, where
Theorem \ref{thm:minimizing-step-function} establishes that minimizing single-well potentials are step functions with explicit control over the number of jumps on any edge.  
Theorem \ref{thm:minimizing-piecewise-linear} then 
uses similar arguments to establish that minimizing convex potentials are piecewise linear, likewise with an explicit estimate on the number of non-smoothness per edge..

In \S \ref{sec:non-constant}
we show that in general the minimizing convex potentials are not constant, in contrast to Lavine's result for the interval.  A specific criterion for this is articulated in Lemma \ref{lem:small}, and it is used to establish specific circumstances when the constant potential is not the minimizer,
in Corollaries \ref{cor:vanishing} and \ref{cor:multiple} and in Theorem \ref{thm:constant-perturbation}.

Finally, in \S \ref{sec:bounds} we show 
through examples
that it is not generally possible to find nontrivial bounds on the fundamental gap in terms of the diameter of the graph alone,
within the given classes.
Example \ref{ex:star-graph-no-bound-1} in 
\S \ref{sec:bounds}
provides a sequence of convex (indeed, linear) potentials on a star graph of fixed diameter for which the fundamental gap $\Gamma$ tends to $0$, while Example \ref{ex:star-graph-no-bound-2} exhibits an analogous sequence of single-well potentials on a star graph.  Conversely,
Example \ref{ex:star-graph-no-bound-3} consists of a sequence of tree graphs of fixed diameter for which $\Gamma$ is unbounded.

\section{Preliminaries: Schr\"odinger operators on metric graphs}\label{sec:Prelims}

We start with the basic framework we will use throughout the paper.

\subsection{Notation}
Given a finite set $\mathcal{E}$ of edges $e$, where each edge $e \in \mathcal{E}$ is identified with a compact interval $[0,\vert e \vert]$ of length $\vert e\vert>0$, whose ends are partitioned into a finite set of equivalence classes, or vertices $v \in \mathcal{V}$, we define the corresponding metric graph $G$ as the natural metric space arising as the union of the edges subject to the adjacency relations imposed by the vertices as described in \cite{Mugnolo}, see also \cite[Section~1.3]{Berkolaiko}. We always assume unless explicitly stated otherwise that $G$ is connected; thus it is a compact, connected metric space which necessarily has finite total length, which we denote $L := \sum_{e\in\mathcal{E}} \vert e \vert$, and finite diameter $D$.
 
We introduce two important functional spaces on $G$, the Lebesgue space $L^{2}(G)$ and the Sobolev space $W^{1,2}(G)$:
\begin{defi}

\noindent We say that $u\in L^{2}(G)$ if $u_{e}\in L^{2}(0,l_{e})$ for all $e\in \mathcal{E}$, where $u$ is a function on metric graph $G$ and $u_{e}$ ia a collection function defined on $(0,l_{e})$ for all $e\in\mathcal{E} $, and 
$$
\Vert u \Vert^{2}_{L^{2}(G)}=\sum_{e\in\mathcal{E}}\Vert u_{e} \Vert^{2}_{L^{2}(0,l_{e})}<\infty .
$$
\end{defi}

\begin{defi}
The Sobolev space $W^{1,2}(G)$ consists of all continuous functions
on $G$ that belong to $W^{1,2}(e)$ for each edge $e$ and such that 
$$
\Vert u \Vert^{2}_{W^{1,2}(G)}=\sum_{e\in\mathcal{E}}\Vert u_{e} \Vert^{2}_{L^{2}(0,l_{e})}+\Vert u'_{e} \Vert^{2}_{L^{2}(0,l_{e})}<\infty.
$$
\end{defi}

We consider on each edge of the graph $G$ the self-adjoint Schr\"odinger operators $H$ defined on an appropriate domain of definition in $L^{2}(G) $, whereby for $u\in W^{1,2}(G)$
$$
Hu=-u''(x)+q(x)u(x).
$$

At each vertex $v \in \mathcal{V}$, we could
impose either a Dirichlet condition $u(v)=0$ or a $\delta$-coupling (or Robin) condition
$$
 \sum_{x_{i}\in v}\partial  u(x_{i})  =\alpha_{v}u(v),
$$
where $\partial u(x_{i})$ is the normal derivative of $u$ on the edge $e_i$ pointing into $v$. However, to keep the exposition simple, by default we will always restrict to the case $\alpha_v = 0$ of Kirchhoff conditions, i.e., $
\displaystyle \sum_{x_{i}\in v}\partial  u(x_{i})  =0,$ unless explicitly stated otherwise, and in particular all our results will be formulated only for this case.

It is well known that the spectrum of this operator is discrete on any compact metric graph \cite[Theorem~3.1.1]{Berkolaiko}. Moreover, we have the following Perron--Frobenius-type result, see \cite{Kurasov}.

\begin{prop}\label{Proposition3-1}
Let $G$ be a finite compact connected metric graph. Then, the first eigenvalue of $H$ is simple and the corresponding eigenfunction may be chosen strictly positive.
\end{prop}

Given a compact graph $G$, we will write
\begin{displaymath}
    \lambda_1(G,q) < \lambda_2(G,q) \leq \ldots
\end{displaymath}
for the eigenvalues, suppressing one or more of the arguments (in particular the graph $G$) as convenient if there is no danger of confusion, and $u_n \sim \lambda_n$, $n\geq 1$, for the corresponding eigenfunctions, chosen to form an orthonormal basis of $L^2(G)$ (thus, in particular, $\|u_n\|_2 = 1$ for all $n$) unless otherwise specified.

It is clear that the spectrum of the Schr\"odinger operator $-\Delta+q$ on a quantum graph $G$ is positive when $q$ is positive and the fundamental gap is unaffected by adding a constant, that is, $\Gamma[q]=\Gamma[q+c]$ for all $c \in \R$, where we recall $\Gamma[q] = \lambda_2 (q) - \lambda_1 (q)$.

We will mostly consider an important subclass of graphs, the class of trees; we recall that a metric graph $G$ is a \emph{tree graph} (or \emph{tree} for short) if $G$ is connected and has no cycles.

\begin{defi}
\label{def:leaf}
We call a vertex $v \in \mathcal{V}(G)$ a \emph{leaf} if $\deg v = 1$. If $G$ is a tree, and only in this case, we call the set of all leaves of $G$ the \emph{boundary} of $G$, denoted by $\partial G$.
\end{defi}

(See \cite[Introduction and Section 4.2]{KeRo} for a discussion of why the set of leaves is often not a good notion of the ``boundary'' of a non-tree quantum graph.)

\subsection{Convex and single-well potentials on metric trees}
\label{sec:convex-single-well}

An initial question to confront is how to interpret the notion of convexity of a function on a metric graph $G$.   

Del Pezzo, Frevenza and Rossi \cite {Rossi2} have considered the question of defining convex functions on metric graphs, proposing that 
$u:G\rightarrow\mathbb{R}$ is convex on $G$ provided that, for any $x,y \in G$, $x \neq y$, it satisfies the condition
\begin{equation}
\label{eq:convex-function-inequality}
u(z)\leqslant \frac{\dist(y,z)}{\dist(x,y)}u(x)+\frac{\dist(x,z)}{\dist(x,y)}u(y),
\end{equation}
by way of analogy with a common characterization of them on intervals and domains.  This definition is equivalent to requiring that the function $u$ be convex on every path $P \subset G$.  
We prefer to impose convexity on only certain paths in $G$, which leaves us with a richer family of spectral problems, including nontrivial functions on cycles, to which our methods apply.

\begin{defi}
\label{def:path}
Let $G$ be a metric graph. 
A \emph{path} $ P \subset G$ is the image in $G$ of a continuous map $\gamma: [0,1] \to G$ which is injective except that $\gamma (0) = \gamma (1)$ is allowed, in which case we call $P$ a \emph{closed path}.
\end{defi}

\begin{defi}\label{definition2-12}
Let $\mathcal{P}$ denote a distinguished finite set of paths $P \subset G$. We say the function $u:G\rightarrow\mathbb{R}$ is \emph{convex with respect to} $\mathcal{P}$ if $u$ is convex on each $P$, regarded as an interval.  
We say $u$ is \emph{strictly convex with respect to} $\mathcal{P}$ if the restriction of $u$ to each $P$ is always strictly convex. 
\end{defi}

If, given any pair of points $x,y \in G$, each path connecting $x$ and $y$ is contained in one of the paths 
$P \in \mathcal{P}$, then convexity with respect to $\mathcal{P}$ reduces to convexity on $\Omega$ in the sense 
of \cite {Rossi2}. In particular:
\begin{lem}
\label{lem:convex-characterization}
Let G be a metric graph and let any $x,y\in G$ with $x\neq
y$. The following are equivalent:
\begin{enumerate}
\item $u$ is convex on each path $P \in \mathcal{P}$.
\item $u$ satisfies the condition (\ref{eq:convex-function-inequality}).

\end{enumerate}
\end{lem}

\begin{rem}
\label{rem:affine}
Given a compact connected metric graph $G$, it is natural to say that $u: G \to \R$ is \emph{affine} if, for all $x,y \in G$ and $z \in [x,y]$ (with $z\neq x,y$),
$$
u(z) = \frac{\dist(y,z)}{\dist(x,y)}u(x)+\frac{\dist(x,z)}{\dist(x,y)}u(y).
$$
It is not hard to see that if $G$ is not a path graph, then the \emph{only} affine functions on $G$ are the constant functions.
\end{rem}

\subsection{Classes of potentials}

With the above background, and as mentioned in the introduction, we can now introduce the two natural and common classes of potentials we will be considering, namely convex and single-well. In order to recover compactness properties of the classes in a suitable sense, and thus existence of optimizers (which will be the subject of Section~\ref{sec:existence}), we will always need to impose a uniform $L^\infty$-bound on the class.

We start with the convex case; here, as discussed above, we may consider any compact graph $G$, not necessarily a tree, and restrict to a suitable subset.

\begin{defi}\label{eq:convex-class-m}
For a given family of paths $\mathcal{P}$ and a given $M>0$ we define
\begin{equation}
\label{bdd-convex-pathwise}
    C_{G,\mathcal{P},M} := \left\{ q: G \to \R: q \text{ is convex on } \mathcal{P} \text{ and }
    0\leq q(x)\leq M \text{ for all } x\in\cup_\mathcal{P}\mathcal{P}\right\}.
\end{equation}
\end{defi}

This set is isomorphic to a subset of the direct sum of the set of convex functions on compact intervals, for each of which the $L^\infty$-norm is bounded by $M$.

We will often consider tree graphs, for which it is natural to choose $\mathcal{P}$ as the set of all paths from one boundary vertex to an other.  In this case we denote $C_{G,\mathcal{P},M}$ simply by $C_{G,M}$. 

The set $C_{G,\mathcal{P},M}$ is far from the only interesting set of potentials commonly considered. Another common category is the ``single-well'' potentials.  Let us recall the definition on an interval:

\begin{defi}
\label{def:single-well-interval}
The function $q$ is called a 
{\rm single-well function} on an interval $[0,l]$ if $q$ is non-increasing on $[0,a]$ and non-decreasing on $[a,l]$,
for some $ a \in[0,l]$.
Any such point $ a $ is called a {\rm transition point} (which need not be unique). 
\end{defi}

Inspired by Lemma~\ref{lem:convex-characterization}, we define single-well functions on compact metric graphs
as follows.

\begin{defi}\label{eq:single-well-class-m}
Given a compact metric tree $T$, we say the function $u:T\rightarrow\mathbb{R}$ is \emph{single-well on $T$} if there exists a point
$a \in T$ such that on every path from a boundary vertex to $a$, $u$ is a nonincreasing monotonic function.  If $G$ is a compact metric graph containing the tree $T$, we say $u:G\rightarrow\mathbb{R}$ is \emph{single-well with respect to the subtree $T$} if it is a single-well function when restricted to $T$.  The set of such functions will be denoted $SW_{G,T,M}$.  If $G=T$, we write $SW_{G,M}$.
\end{defi}

We identify single-well functions which agree almost everywhere. One can show that $C_{G,M} \subset SW_{G,M}$ for all $M > 0$ if $G$ is a tree, but we will not need this. It is clear that $SW_{G,T,M_1} \subset SW_{G,T,M_2}$ and $C_{G,\mathcal{P},M_1} \subset C_{G,\mathcal{P},M_2}$ for all $M_2 \geq M_1 > 0$.

One could define families of $N$-well functions on graphs in a natural analogous way.

\subsection{Perturbation of eigenvalues} 
\label{sec:feynman-hellman}

Before proceeding to discuss our main results, for future reference we collect a few important results on the behavior of the spectral gap with respect to perturbations of the potential, for a fixed graph. The first is a simple continuity result; the second is a variation of an explicit formula for the first derivative of a simple eigenvalue with respect to perturbations of $q$, often referred to as a Feynman--Hellman formula. This played a key role in the works of Lavine \cite{Lavine} and other authors in the case of intervals, and extends immediately to quantum graphs.

However, on graphs the second eigenvalue need not be simple, thus we need a modified version which allows for degeneracies. 
Although degenerate perturbation theory is standard lore in
quantum physics, and
addressed carefully in the works of Kato (e.g., \cite{Kato} VII \S4), for reference
we formulate here a proposition 
that can be readily applied to the fundamental gap $\Gamma = \lambda_2 - \lambda_1$.

Notationally, we fix a compact metric graph $G$ and a suitable potential energy $q$ on $G$.

\begin{lem}
\label{lem:potential-continuity}
If $G$ is compact, the eigenvalues of $H$ depend continuously on $q$ in the $L^p(G)$ topology for all 
$1 \le p \le \infty$.
\end{lem}
\begin{proof}
According to a theorem of Kato, to establish type (B) analytic behavior it suffices to have an estimate of the form
\begin{equation}
\left|\int_G q(x) |\zeta(x)|^2\,{\rm d}s\right| \le a \int_G |\zeta'(x)|^2\,{\rm d}s + b  \int_G |\zeta(x)|^2\,{\rm d}s
\end{equation}
for functions $\zeta \in W^{1,2}(G)$ and any finite $a,b$ \cite{Kato,RSIV}.  This is obvious for $q \in L^\infty(G)$ and by interpolation will be true for all $q \in L^p(G)$ if it is established for $p=1$.  Since $W^{1,2}(G)$ is isomorphic to a closed subset of $H^1(I)$ for some interval $I$, it inherits the Sobolev property of being compactly embedded in $L^\infty (G)$.  Therefore 
\begin{equation}
\left|\int_G q(x) |\zeta(x)|^2\,{\rm d}s\right| \le \|q\|_{L^1(G)} \|\zeta\|_{L^\infty(G)}^2  \le C \|q\|_{L^1(G)}  \|\zeta\|_{W^{1,2}(G)}^2,
\end{equation}
as required.
\end{proof}

For the explicit perturbation formula, we first introduce some terminology. A perturbation $P$ will be called 
{\em positively admissible} if for sufficiently small $t > 0$, $q+t P$
remains within the feasible set for the optimization problem (for example, convex or single-well
potentials).  It will be termed {\em admissible} if 
both $P$ and $-P$ are positively admissible.  We will denote by $\{u_2^{(k)}\}$ an orthonormal basis for the eigenspace of $\lambda_2(q)$, and 
as usual will let $u_1 = u_1(q)$ denote an orthonormalized eigenfunction for $\lambda_1$.  To simplify notation we assume that 
$u_1$,  $u_2^{(k)}$, and $P$ are real-valued.

\begin{prop}\label{degpert}
Suppose that $q_*$ is a minimizing potential for the function $\lambda_2 (q) - \lambda_1 (q)$ within a given class $\mathcal{C} \subset L^\infty (G)$, and, supposing $\lambda_2(q_*)$ is an eigenvalue with
finite multiplicity $n \geq 1$, for a given function $P$ on $G$ let the $n \times n$-matrix $(P_{jk})$ be defined by
$$P_{jk} := \int_G{P(x) \left(u_2^{(j)}(x) u_2^{(k)}(x) - (u_1(x))^2\right)}\,\dx .$$
\begin{enumerate} 
\item
If $P$ is a bounded positively admissible perturbation, then $P_{jk}$ is a non-negative matrix.  In particular, for each $k$,
$P_{kk} \ge 0$.
\item
If $P$ is a bounded admissible perturbation, then $P_{jk}$ is the 0 matrix.
\end{enumerate} 
Suppose that $q_*$ is a maximizing potential within the class $\mathcal{C}$. 
\begin{enumerate} 
\item[(3)]
If $P$ is a bounded positively admissible perturbation, then it is not possible for $P_{jk}$ to be a strictly positive matrix, i,e, at least one of its eigenvalues must be $\le 0$.
\item[(4)]
If $P$ is a bounded admissible perturbation, then it is not possible for either $P_{jk}$ or $-P_{jk}$ to be a strictly positive matrix.
\end{enumerate} 

\end{prop}

\begin{rem}
This proposition is easy to adapt to other functions $F(\lambda_1, \lambda_2)$ defined on the spectrum and to other sorts of self-adjoint operators.  For example, 
it applies to $\lambda_2 - \alpha \lambda_1$ for any constant $\alpha$, notably $\alpha=0$.
\end{rem}

\begin{proof}
Since this is standard lore in quantum mechanics textbooks, e.g., \cite{Schiff}, we content ourselves with a sketch.  According to regular perturbation theory
\cite{Kato}, there are eigenvalues $\lambda_2^{(k)}(t)$ corresponding
to the potential $q_* + t P$ which are analytic in $t$ near $0$ and coincide with $\lambda_2(q_*)$ when $t=0$.  The derivatives 
$(\lambda_2^{(k)})'(0)$ are the eigenvalues of $P_{jk}$.  Since by definition $\lambda_2(t) = \min_k \lambda_2^{(k)}(t)$, the minimality of $\Gamma(q_*)$ would be contradicted if even one of these derivatives is 
respectively negative in case (1) or nonzero in case (2).  

In the maximizing cases, all of the derivatives $(\lambda_2^{(k)})'(0)$ would have to be $>0$ to obtain a contradiction, hence the stronger condition on 
the matrix $P_{jk}$.
\end{proof}

We formulate a simple corollary of Proposition~\ref{degpert} which gives a useful practical criterion for showing non-optimality of a given potential $q$ within some class $\mathcal{C}$, without requiring $\lambda_2 (q)$ to be simple.

\begin{cor}
\label{cor:feynman-hellman}
Let $G$ be a compact connected graph, and let $\mathcal{C} \subset L^\infty (G)$ be a class of potentials. Given a potential $q \in \mathcal{C}$, suppose that $P$ is a perturbation such that $q+tP \in \mathcal{C}$ for each $t \geq 0$, respectively $t \leq 0$, sufficiently small. If there exist normalized eigenfunctions $u_1$ associated with $\lambda_1(q)$ and $u_2$ associated with $\lambda_2(q)$ such that
\begin{equation}
\label{eq:feynman-hellman-condition}
    \int_{G}P(x)[u_{2}^{2}(x)-u_{1}^{2}(x)]\, \dx < 0,
    \text{ respectively, } > 0,
\end{equation}
then $q$ does not minimize the fundamental gap in $\mathcal{C}$.
\end{cor}

\section{Compact classes of potentials and existence of optimizers}
\label{sec:existence}

\subsection{Compactness properties of convex and single-well classes}
\label{sec:compactness}

In this section we will show that within the classes $C_{G,\mathcal{P},M}$ and $SW_{G,T,M}$ of convex and single-well potentials (for a given compact graph $G$, given set of paths $\mathcal{P} \subset G$, respectively a given subtree $T \subset G$, and given $M>0$) there are always potentials minimizing and maximizing the spectral gap.

This is based the fact that each class enjoys a suitable compactness property; in each case the underlying heuristic principle is the same, but details are slightly different, using either the Blaschke selection principle or the Helly compactness theorem. In particular, it would be possible to prove similar compactness properties for other classes of potentials (such as, e.g., $N$-well functions), but the proof in each case would need to be slightly adjusted.

\begin{lem}\label{lemma2-10}
Given $G$, $\mathcal{P} \subset G$ and $M>0$, the set $C_{G,\mathcal{P},M}$ is sequentially compact in the uniform topology. 
\end{lem}

\begin{proof}
This is a direct consequence of the Blaschke selection principle, which implies that a uniformly bounded sequence of convex functions on a compact interval contains a subsequence that converges uniformly to a limit, which is convex. Given that $\mathcal{P}$ may be identified with a finite union of intervals (i.e. the paths), and any $q \in C_{G,\mathcal{P},M}$ is convex when restricted to each path, We may apply the principle sequentially on each, passing to a further subsequence in each case. The resulting limit function will be convex on every path, and thus on $\mathcal{P}$.
\end{proof}

In the context of single-well potentials, we recall that according to the Helly compactness theorem, a uniformly bounded sequence of monotonic functions on an interval has a subsequence that converges pointwise to a monotone limit.  It is not hard to extend this convergence property to uniformly bounded single-well functions on an interval, which are therefore compact in the $L^p$ topology for every $p < \infty$ \cite{ElAHar}.
In the context of metric trees (or subtrees of a graph), a similar compactness can be derived from Helly's theorem, which can be applied due to the monotonicity properties of single-well functions along each path emanating from the bottom of the ``well.'' In any uniformly bounded sequence of single-well functions, it is possible to extract a convergent subsequence that maintains the single-well property.

\begin{lem}
\label{lem:sw-compact}
    Given $G$, $T \subset G$ and $M>0$, the set $SW_{G,T,M}$ is sequentially compact in $L^p(G)$ for any $1 \leq p < \infty$.
\end{lem}

\begin{proof}
The argument is identical to the proof of Lemma~\ref{lemma2-10}, except that, as mentioned, one replaces the Blaschke selection principle with the Helly compactness theorem in the form of \cite[Proposition~2.2]{ElAHar}.
\end{proof}

\subsection{Existence of optimizers}

We can now establish the existence of potentials minimizing and maximizing the spectral gap within classes of the form
$C_{G,\mathcal{P},M}$ and $SW_{G,T,M}$. We stress that the method, which basically combines the continuity result of Lemma~\ref{lem:potential-continuity} with the compactness results of Section~\ref{sec:compactness}, 
is quite general and could be easily adapted to prove existence of optimizers in other classes such as, for instance, $N$-well potentials. Likewise, we will focus on the fundamental gap, but the same method could also be applied to other interesting functions of the spectrum such as sums of the first $k$ eigenvalues, or partition functions.

\begin{prop}\label{proposition3-2}
Let $\mathcal{F} \subset L^p(G)$ be compact in the $\|\cdot \|_{L^p}$ topology, and let $f(\lambda_1, ..., \lambda_{\ell})$ be a differentiable function of the first $\ell$ eigenvalues of $H$, $\ell < \infty$.  Then $f$ achieves its maximum and minimum on $\mathcal{F}$, and the first variations of $f$ with respect to perturbations in $\mathcal{F}$ can be calculated by the procedures of formal perturbation theory.
\end{prop}

\begin{proof}
This follows directly from Lemma~\ref{lem:potential-continuity} and standard theory \cite{Kato}.
\end{proof}

\begin{prop}\label{existenceprop}
Let $G$ be a compact metric
graph equipped with standard vertex conditions and fix $1 \le p < \infty$ and $M > 0$. Then $C_{G,\mathcal{P},M}$ and $SW_{G,T,M}$ are compact subsets of $L^p(G)$.
\end{prop}

\begin{proof}  
Since $G$ is compact, $L^\infty (G) \subset L^p(G)$, and the proof follows immediately from Lemmas~\ref{lemma2-10} and~\ref{lem:sw-compact}, respectively.
\end{proof}

\begin{thm}{\label{thm:exist}}
Let $G$ be a compact 
metric graph equipped with standard vertex conditions, fix a family of paths $\mathcal{P} \subset G$, a tree subgraph $T \subset G$, and $M>0$. Let $\Lambda$ be either $C_{G,\mathcal{P},M}$ or $SW_{G,T,M}$. Then there exist potentials $q_\ast,q^\ast \in \Lambda$ such that
\begin{displaymath}
    \Gamma[q_\ast] = \min \{\Gamma[q]: q \in \Lambda\}, \qquad
    \Gamma[q^\ast] = \max \{\Gamma[q]: q \in \Lambda\}.
\end{displaymath}
\end{thm}

\begin{proof}
The proof is an immediate consequence of the compactness result Proposition~\ref{existenceprop} and the continuity/smoothness result Proposition~\ref{proposition3-2} applied to $f(\lambda_1,\lambda_2) = \lambda_2 - \lambda_1$.
\end{proof}

\section{Piecewise linearity of minimizing potentials}
\label{sec:characterization}

In this section, we look at the explicit form of the gap-minimizing potential of $H$ as a complement to the existence results of the previous section; here we will generally follow the strategy of \cite{ElAHar}. However, here the conclusion is necessarily much weaker: any minimizing potential must be piecewise linear in the convex case, or piecewise constant in the single-well case. We will see in the following sections that the minimizer need not, however, be constant.

The proofs in the two cases of single-well and convex potentials are structurally quite similar; namely, one uses well-chosen perturbations to exclude zones of strict convexity in the convex case, or anything other than jumps at certain distinguished points of the graph in the single-well case. However, despite these similarities, there does not seem to be a clear way to derive the two results (Theorems~\ref{thm:minimizing-step-function} and~\ref{thm:minimizing-piecewise-linear}) from a single, more general principle. Thus we will give a relatively detailed proof in the single-well case, and then sketch the convex case.

In both cases we will need the following technical result.

\begin{lem}{\label{Lemma4.3}}
Let $G$ be a compact metric graph and let $q\in L^{1}_{loc}(G)$. Suppose $u_1$ and $u_k$ are, respectively, normalized eigenfunctions of $\lambda_{1} = \lambda_1 (q)$ and $\lambda_{k} = \lambda_k (q)$. Then, for any edge $e \in \mathcal{E}$:
\begin{enumerate}
    \item The closure of any connected component of $e \cap \{u_k \neq 0\}$ contains at most two zeros of $u_{k}^{2}-u_{1}^{2}$; in particular, if $u_k$ does not change sign in the interior of $e$ then there are at most two zeros in $e$ in total;
    \item  If $v \in \mathcal{V}$ is a leaf of $G$, $e \sim v$, and $u_k \not\equiv 0$ on $e$, then the connected component of $e \cap \{u_k \neq 0\}$ which contains $v$, contains at most one zero of $u_{k}^{2}-u_{1}^{2}$.
\end{enumerate}
In particular, the zeros of $u_{k}^{2}-u_{1}^{2}$ are isolated.
\end{lem}

Note that edges are taken to be closed; that is, a zero of $u_k^2-u_1^2$ at a vertex counts towards the total on any adjacent edge.

\begin{rems}
(1) Lemma~\ref{Lemma4.3} is sharp in the base case where $G$ is just an interval and $q=0$: in this case $u_1$ is constant and $u_k$ is sinusoidal; the normalization guarantees that on each nodal domain of $u_k$, $\max |u_k| > |u_1|$ and so $u_k^2-u_1^2$ must have exactly two zeros on that nodal domain \emph{unless} the nodal domain is adjacent to one of the endpoints of the interval, in which case the Neumann condition guarantees there will be exactly one. We note in particular that the proof of (1) does not use the Kirchhoff vertex conditions that $u_1$ and $u_k$ satisfy.

(2) In the case that $k=2$, $G$ is a tree, $q$ belongs to one of the classes $C_{G,M}$ or $SW_{G,M}$, and standard conditions are assumed, it seems quite possible that the number of zeros of $u_2^2 - u_1^2$ could nevertheless be controlled much more closely. Essentially, the function $\frac{u_2}{u_1}$ (studied in the proof of the lemma) should be monotonic on each path leading away from the zero set of $u_2$, as a closer study of the Wronskian $W$ associated with $u_1$ and $u_2$ plus the behavior of $\frac{u_2}{u_1}$ across the vertices should show. We will not attempt to improve the estimate here, though.
\end{rems}

\begin{proof}[Proof of Lemma~\ref{Lemma4.3}]
We adapt the approach of \cite[Lemma~3.2]{ElAHar}, in turn inspired by \cite{Ash1},  although only in case (2) can their approach be applied directly.

(2) We start with the case more similar to the interval. We assume without loss of generality that $u_k(v)>0$, and take the closure of $e \cap \{u_k \neq 0\}$ to be parametrized as an interval of the form $[0,\ell]$, where $0$ corresponds to $v$. We suppose for a contradiction that there exist $\xi_1,\xi_2 \in [0,\ell]$, $\xi_1 < \xi_2$, such that
\begin{displaymath}
    |u_k(\xi_i)| = |u_1(\xi_i)| = u_1 (\xi_i)
\end{displaymath}
for $i=1,2$. We may define a new function $w:=\frac{u_k}{u_1} > 0$, so that $w(\xi_1) = w(\xi_2) = 1$. By Rolle's theorem, there exists $\zeta \in (\xi_1,\xi_2)$ such that $w'(\zeta) = 0$. On the other hand, denoting by
\begin{displaymath}
    W(x) = u_1(x)u_k'(x) - u_k(x)u_1'(x),
\end{displaymath}
the Wronskian associated with $u_1$ and $u_k$ in $x \in (0,\ell)$, we find using the eigenvalue equation that
\begin{displaymath}
    W'(x) = (\lambda_1 - \lambda_k)u_1(x)u_k(x) < 0
\end{displaymath}
for all $x \in (0,\ell)$ (by assumption on $u_k$ and the fact that $u_1 > 0$ on $G$), and thus, since $W(0) = 0$ by the Kirchhoff condition that $u_1$ and $u_k$ satisfy at $v$,
\begin{displaymath}
    w'(x) = \frac{W(x)}{u_1^2(x)} = \frac{\lambda_1 - \lambda_k}{u_1^2(x)}
    \int_0^x u_1(t)u_k(t)\,\dt < 0
\end{displaymath}
for all $x \in (0,\ell)$, a contradiction.

(1) Suppose now that $e$ is any edge on which $u_k$ does not vanish identically (if it is identically zero, then since $u_1$ does not vanish anywhere, $u_k^2-u_1^2$ cannot vanish at all on $e$), and let $I \subset e$ is a closed interval on which $u_k$ does not change sign, although it is allowed to vanish at the endpoints of $I$. We may assume without loss of generality that $u_k > 0$ in the interior of $I$.
For an arbitrary but fixed orientation of $e$ and hence of $I$, define $w$ and $W$ in the interior of $I$ as in (2), and observe that $w'(x) = \frac{W(x)}{u_1^2(x)}$ still, so that the zeros of $w'$ and $W$ coincide. Now $W$ can have at most one zero in $I$. Indeed, the eigenvalue equation still implies that
\begin{displaymath}
    W'(x) = (\lambda_1 - \lambda_k)u_1(x)u_k(x)
\end{displaymath}
in the interior of $I$, which is never equal to zero there by assumption on $u_k$. It follows that $W$ is strictly monotonic on $I$ and hence has at most one zero. But then $w'$ likewise can have at most one zero in the interior of $I$.

Now if there were three points $\xi_i \in I$, $i=1,2,3$, such that $u_k(\xi_i) = u_1(\xi_i)$, with $e$ parametrized in such a way that $\xi_1 < \xi_2 < \xi_3$, then by Rolle's theorem $w'$ would have at least one zero in $(\xi_1,\xi_2)$ and at least one zero in $(\xi_2,\xi_3)$, a contradiction.
\end{proof}

\subsection{Single-well potentials}
We first deal with minimizers in classes of single-well potentials, whose existence is guaranteed by Theorem~\ref{thm:exist}.

\begin{thm}
\label{thm:minimizing-step-function}
Let $T$ be a tree and $M>0$. If $q_\ast$ minimizes $\Gamma[q]$ in the set $SW_{T,M}$, then $q_\ast$ is a non-constant step function. Moreover, $\lambda_2 (q_\ast)$ is simple, and, taking $u_1$ and $u_2$ as normalized eigenfunctions of $\lambda_1 (q_\ast)$ and $\lambda_2 (q_\ast)$, respectively, $q_\ast$ can jump only at a zero of $u_2^2-u_1^2$. In particular, assuming without loss of generality that $u_2$ changes sign only at a (necessarily unique) vertex, on any given edge, $q_\ast$ has at most two points of discontinuity in the interior of that edge, or at most one if the edge ends at a leaf.
\end{thm}

The theorem holds for general graphs $G$ and $SW_{G,T,M}$ as long as $T$ is a \emph{spanning tree} of $G$.
It seems unlikely that the assertion that $q_\ast$ has at most two discontinuities per edge (which comes from Lemma~\ref{Lemma4.3}) is optimal; we expect it can be at most one, namely jumping once between $0$ and $M$.

\begin{proof}
Denote by $x_m \in T$ any point at which $q_\ast (x_m) = \displaystyle\min_{x \in T} q_\ast (x)$, which we may assume to be a vertex (since any point in the interior of an edge can be designated as a degree two vertex in the usual way). Then $q_\ast$ is non-decreasing on any path from $x_m$ to any leaf of $T$.

We first show that $q_\ast$ cannot be constant. If it were, then after a shift we may assume $q_\ast = \frac{M}{2}$ on $T$. Let $I \subset T$ be any interval (connected subset of an edge) of positive length such that $u_{2}^{2}-u_{1}^{2}>0$ on $I$; the existence of such an interval is guaranteed by Lemma~\ref{Lemma4.3} and the fact that $u_1$ and $u_2$ both have $L^2$-norm $1$. We define a family of potentials
\begin{displaymath}
    q_t := q_\ast + t\chi_I,
\end{displaymath}
then $q_t$ is single-well for $t\leq 0$, and still in $SW_{T,M}$ for $t\in [-\frac{M}{2},0]$. Since
\begin{displaymath}
    \int_{I} u_{2}^{2}-u_{1}^{2}\, \dx >0,
\end{displaymath}
by Corollary~\ref{cor:feynman-hellman} we immediately obtain that $\Gamma[q_t] < \Gamma[q_\ast]$ for $t<0$ small, a contradiction to minimality of $q_\ast$.

We will next give the argument that, for any normalized second eigenfunction $u_2$, $q_\ast$ is constant on \emph{any} connected subset of $T$ on which $u_{2}^{2}-u_{1}^{2}$ does not change sign.

Indeed, suppose $S \subset T$ is a closed, connected set with, say, $u_{2}^{2}-u_{1}^{2}>0$ on $S$, and suppose for a contradiction that $q_\ast$ is not constant on $S$. We may assume without loss of generality that there exists a point $0 \in S$ which is either a leaf of $T$ or a boundary point of $S$ in $T$, such that $q_\ast$ is monotonically increasing on every path in $S$ starting at $0$. Let $a$ be any point at which $q_\ast$ attains its maximum in $S$, necessarily also a leaf or a boundary point.

For $t \in [0,1]$ we choose a perturbing family of the form
\begin{displaymath}
    q_t (x) := tq_\ast (0) + (1-t)q_\ast (x)
\end{displaymath}
for all $x \in S$, and $q_t =q_\ast$ on the rest of $T$. Then $q_t$ is still single-well on $T$ (being clearly non-decreasing on every path in $S$ emanating from $0$, with $q_t(0)=q_\ast(0)<q_t(x)\leq q_\ast (x)$); but by Corollary~\ref{cor:feynman-hellman},
\begin{displaymath}
    \Gamma_+'[q_t]|_{t=0} \leq \int_{S}(u_{2}^{2}-u_{1}^{2})(q_\ast(0)-q_{\ast}(x))\,\dx \leq 0
\end{displaymath}
(where the derivative is with respect to $t$). Our assumption on $u_1$ and $u_2$ plus the fact that $q_\ast$ is minimizing, and all functions are continuous, implies that $q_\ast (x) = q_\ast (0)$ for all $x \in S$, as claimed.

If $u_{2}^{2}-u_{1}^{2}<0$ on $S$, then the argument is slightly more complicated. Noting that $S$ is itself a tree, we take $a \in S$ to be any point in $\{x \in S: q_\ast (x) = \max_S q_\ast (x)\}$ closest to $0$, and observe that, necessarily, regardless of whether $a$ is a vertex, there exists a unique path in $S$ from $0$ to $a$, and an interval (subset of an edge) $I \subset S$ on that path which ends at $a$. In particular, $q_\ast(x) < q_\ast (a)$ for all $x \in I \setminus \{a\}$. We may then consider $q_t(x) = tq_\ast (a) + (1-t)q_\ast (x)$ on $I$, and $q_t = q_\ast$ on $T \setminus I$. Running through the above argument, suitably adjusted, leads to $\Gamma_-'[q_t]|_{t=0} \geq 0$, and thus $q_\ast(x) = q_\ast (a)$ on $I$, a contradiction to our choice of $a$. We conclude that $q_\ast$ is constant on $S$.

We next show that $\lambda_2 (T,q_\ast)$ is simple. Since $q_\ast$ is non-constant, it must have at least one jump, say at $x_\ast \in T$, so that $|u_2(x_\ast)| = |u_1(x_\ast)|$ for some eigenfunction $u_2$.

Suppose there exists a linearly independent normalized eigenfunction $\tilde u_2$, then either $\tilde u_2 (x_\ast) = 0$, in which case $|\tilde u_2(x)| < |u_1 (x)|$ in a neighborhood of $x_\ast$ and we can run the above argument to conclude that $q_\ast$ is constant around $x_\ast$, or else $\tilde u_2 (x_\ast) \neq 0$, in which case $(u_2+t\tilde u_2)^2 - u_1^2$ will have a zero close to, but not at, $x_\ast$ for $t \neq 0$ small enough, again meaning $q_\ast$ cannot have a jump at $x_\ast$, a contradiction. Hence $\lambda_2 (T,q_\ast)$ is simple.

To finish the proof, we observe that since $G$ is a tree $u_2$ can change sign at most once on any given edge (indeed, either $u_2$ changes sign at a vertex or there is exactly one edge on which it changes sign, and we declare this point to be a dummy vertex). By Lemma~\ref{Lemma4.3}, $u_2^2-u_1^2$ can be zero at most twice per edge; thus $q_\ast$ can jump at most twice per edge, or once if the edge ends at a degree-one vertex.
\end{proof}

\begin{rem}
\label{rem:single-well-maximizer}
If $q^\ast$ is a maximizing potential for $SW_{G,T,M}$, and in addition $\lambda_2 (T, q^\ast)$ is simple, then the stronger conclusion of Proposition~\ref{degpert} for the minimizers, suitably adjusted, also trivially holds for $q^\ast$, that is,
\begin{displaymath}
    \int_G P(x) \left( u_2(x)^2 - u_1(x)^2\right)\,\dx \leq 0
\end{displaymath}
for every bounded positively admissible perturbation $P$ of $q^\ast$. With this result in hand, the proof of Theorem~\ref{thm:minimizing-step-function} can be directly adapted to show that \emph{if} $\lambda_2 (T, q^\ast)$ is simple, \emph{then} $q_\ast$ is a step function which can only jump at zeros of $u_2^2-u_1^2$. We leave it as an open problem to investigate whether $\lambda_2 (T, q^\ast)$ is necessarily simple if $q^\ast$ is a maximizer.
\end{rem}

\subsection{Convex potentials}
We next deal with convex potentials. We emphasize that here one cannot expect the minimizers to be constant, see Section~\ref{sec:non-constant}.

\begin{thm}
\label{thm:minimizing-piecewise-linear}
If $q_\ast$ minimizes $\Gamma[q] = \lambda_{2}(G,q)-\lambda_{1}(G,q)$ in the category $C_{G,\mathcal{P},M}$, then $q_\ast$ is piecewise linear on $\mathcal{P}$, with at most two points of non-smoothness in the interior of each edge in $\mathcal{P}$.
\end{thm}

(See Remark~\ref{rem:convex-maximizer} (1) on where the points of non-smoothness can occur.)

\begin{proof}
Considering the case of an interval, it was shown in 
\cite{ElAHar} that $q_{\ast}$ cannot be strictly convex on any open interval, as else a perturbation is possible remaining within $C_{G,\mathcal{P},M}$ for which the Feynman--Hellman formula gives a nonzero derivative.

Let $u_1$ and $u_2$ be any normalized eigenfunctions for $\lambda_1 (G,q_\ast)$ and $\lambda_2 (G,q_\ast)$, respectively; we do not assume that $\lambda_2 (G,q_\ast)$ is simple.

If $q_{\ast}$ were strictly convex on any interval of an edge, and $q_{\ast}<M$, then it could be perturbed by a sufficiently small positive or negative multiple of a $C^{2}$ function $f$ of compact support in that interval. From Lemma~\ref{Lemma4.3} we can always find a subinterval such that $u_{2}^{2}-u_{1}^{2}\neq 0.$  Since then
$$
\int_{ \mathcal{P} }(u_{2}^{2}-u_{1}^{2})f(x)\,\dx\neq 0.
$$
Corollary~\ref{cor:feynman-hellman} yields a contradiction to the minimality of $q_\ast$. It follows that $q_\ast$ is not strictly convex on any interval on $\mathcal{P}$.

Since a convex function on an interval has left and right derivatives at every point, with at most a countable number of values of $x$ for which the two derivatives are not equal \cite{Webster}, we can thus conclude that the distributional second derivative 
of $q_{\ast}$ is of the form
\begin{equation}
\label{eq:minimizing-potential-second-derivative}
    q''_{\ast}(x)=\sum_{j}\alpha_{j}\delta(x-x_{j}),
\end{equation}
 where $\alpha_{j}\geqslant0$.

Arguments adapted from \cite[Proof of Proposition~3.4]{ElAHar} 
can now eliminate the possibility that with the minimizing 
potential energy $q_*$ there could be 
other points of non-smoothness besides 
vertices and at most one point in any connected subset of an edge where $u_2^2 - u_1^2>0$:

Suppose first that  $u_2^2 - u_1^2 < 0$
on an open subset of an edge which contains a point of non-smoothness, which we 
designate $x_1$, i.e., Eq. 
\eqref{eq:minimizing-potential-second-derivative} contains a contribution
$\alpha_1 \delta_{x_1}$ with $\alpha_1 > 0$.  
Let $T(x)$ denote the standard ``tent function'' supported in the interval $[-1, 1]$, on which
$T(x) = 1- |x|$.
For sufficiently small $\varepsilon, T\left(\frac{x-x_1}{\varepsilon}\right)$ will be supported in
an interval on which $u_2^2 - u_1^2 < 0$, and for sufficiently small $t>0$,
$q_* + t T\left(\frac{x-x_1}{\varepsilon}\right)$ will be convex.  
Since $T\left(\frac{x-x_1}{\varepsilon}\right) > 0$ for $x_1 - \varepsilon < x < x_1 + \varepsilon$, 
Corollary \ref{cor:feynman-hellman} would be contradicted.

Next suppose that some subinterval of an edge on which 
 $u_2^2 - u_1^2 \ge 0$
contains two points of non-smoothness of $q_*$, 
or one point of nonsmoothness and a terminal vertex.
We
parametrize them such that  $x_1 < x_2$, where $x_2$ is possibly a terminal vertex.   This time we choose to 
perturb $q_* \to q_* - t T\left(\frac{|2x -x_2 - x_1|}{x_2 - x_1} \right)$, 
which preserves convexity for sufficiently small $t>0$,
and again obtain a contradiction to Corollary \ref{cor:feynman-hellman}.
\end{proof}

\begin{rems}
\label{rem:convex-maximizer}
(1) The proof actually shows that points of non-smoothness occur only at vertices or where $u_2^2 - u_1^2 \ge 0$, and on any edge there can be at most one point of non-smoothness in any connected subset where $u_2^2 - u_1^2 \ge 0$. Moreover, on any edge adjacent to a leaf $v$, there are no points of non-smoothness in any a neighborhood of $v$ on which $u_2^2 - u_1^2 \ge 0$. The statement in the theorem follows immediately from this more precise version and Lemma~\ref{Lemma4.3}, which shows (together with the fact that $u_2$ changes sign at exactly one point) that $\{u_2^2 - u_1^2 \ge 0\}$ can have at most two connected components on any edge. However, as with Theorem~\ref{thm:minimizing-step-function} (and indeed Lemma~\ref{Lemma4.3}), it seems entirely possible that our upper bound of two points of non-smoothness per edge is not sharp

(2) As in the single-well case (see Remark~\ref{rem:single-well-maximizer}), if $q^\ast$ is a maximizing potential for $C_{G,\mathcal{P},M}$, then under the additional assumption that $\lambda_2 (G,q^\ast)$ is simple, so that the stronger form of Proposition~\ref{degpert} is available, the above proof can be directly adapted to show that on every path $P \in \mathcal{P}$ the maximizing potential $q_\ast$ is only piecewise linear, with only a finite number of points of non-smoothness. We do not go into details, and also again leave it as an open problem to investigate the simplicity of $\lambda_2 (G,q^\ast)$.
\end{rems}

\section{Non-optimality of constants in the class of convex potentials}
\label{sec:non-constant}

Given that, for any $G$ and any tree $T \subset G$, the minimizing and maximizing potentials in the class $C_{G,T,M}$ should always be piecewise linear, the question naturally arises as to whether these potentials can be characterized more explicitly, and whether a corresponding bound on the (minimal or maximal) fundamental gap can be given.
We recall that on intervals (and convex domains) the constant potential minimizes the fundamental gap among all convex potentials: these are, respectively, the results of Lavine \cite{Lavine}, and Andrews and Clutterbuck \cite{Andrews}.

Lavine's result may be equivalently reformulated for graphs as follows: given any \emph{path graph} $G$ of length $L=L(G)>0$, which may be identified with the interval $[0,L]$, the fundamental gap of a Schr{\"o}dinger operator with convex potential $q$ and standard vertex conditions (a.k.a.\ Neumann boundary conditions on an interval) satisfies
\begin{displaymath}
    \Gamma [q] \geq \Gamma [c] = \frac{\pi^2}{L^2}
\end{displaymath}
for any constant potential $q(x)=c$ for all $x \in G$.

The purpose of this section is to illustrate that, in general, the constant potential does \emph{not} minimize the fundamental gap, even among convex potentials, on a tree. (Here we \emph{always} consider standard vertex conditions unless stated otherwise.)

In fact, we suspect the property of the constant potential not minimizing the fundamental gap is \emph{generic} with respect to the edge lengths of the graph in the usual sense, and will present evidence below to this effect. (See, e.g., \cite{BerkolaikoLiu} or \cite[Definition~2.1]{KeRo} for more on the notion of generic properties in the context of metric graphs.)

\begin{conj}
\label{conj:constant}
Given any underlying discrete tree graph which is not a path graph, the set of edge-length vectors and $M>0$ for which, on the corresponding metric tree $T$, the constant potential minimizes the spectral gap in the set $C_{T,M}$, is of the first Baire category in $\R^{|\mathcal{E}|}_{+} \times \R_{+}$.
\end{conj}

We already know (Theorem~\ref{thm:minimizing-step-function}) that the constant potentials are \emph{never} minimizers in $SW_{T,M}$ (for $M>0$), something already known for the interval (see \cite[Theorem~3.1]{ElAHar}).

We note in passing that in \cite{ASHBAUGH}, Ashbaugh and Kielty make the conjecture that for convex potentials and (non-symmetric) Robin conditions on an interval, the minimizing potential will be linear but not constant; they prove that this is the case for a Dirichlet condition at one endpoint and a Neumann condition at the other. At the same time, there is a certain connection between standard eigenfunctions on a tree and Robin eigenfunctions on an interval, as the former restricted to an interval will necessarily satisfy some Robin condition at the interval endpoints (an observation frequently used in ``graph surgery'', cf.\ \cite[Section~3.1]{Berkolaiko-1}, which we will use in Section~\ref{sec:bounds}). Thus this may be viewed as weak supporting evidence for the conjecture.

Here, we will present two classes of examples for which the constant potential is not optimal; both are based on the Feynman--Hellman formula of Section~\ref{sec:feynman-hellman}, but will involve choosing different kinds of perturbations of the constant potential.

Actually, we know of no examples other than path graphs (i.e. intervals) for which the constant potential is minimal in either of these classes.

We also note that a continuity argument, which we will not perform here, shows that if for a given tree $G$ (and given $M>0$) the constant potential is not minimal for $C_{G,M}$, then there exists some $\varepsilon_0 > 0$ such that for any graph $G'$ with the same topology as $G$ and whose edge lengths differ from those of $G$ by less than $\varepsilon_0$, the constant potential will also not be minimal in $C_{G_\varepsilon,M}$.

\subsection{Eigenfunctions which are small at a leaf}
We start with the principle that if there is a leaf at which some eigenfunction $u_2$ associated with $\lambda_2(G,0)$ is ``small'' in absolute value at that leaf, then the constant potential is not optimal.

\begin{lem}
\label{lem:small}
Let $G$ be a compact tree of total length $L>0$. If there exists an eigenfunction $u_2$ associated with the Laplacian $(q=0)$ on $G$ and a vertex $v \in \mathcal{V}(G)$ with $\deg(v)=1$ such that
\begin{equation}
\label{eq:small}
    u_{2}^{2}(v)<\frac{1}{L}\int_{G}u_{2}^{2}(x)\,{\rm d}x,
\end{equation}
then, for any $M>0$, the potential $q=0$ is not a minimizer of $\Gamma [q]$ in $C_{G,M}$. In particular, this conclusion holds if $u_{2}(v)=0$ for some vertex $v$ with $\deg(v)=1$.
\end{lem}

(We recall that standard vertex conditions 
are assumed throughout; thus in practice $u_2(v)=0$ at a leaf $v$ implies $u_2$ vanishes identically on the corresponding edge.)

\begin{proof}[Proof of Lemma~\ref{lem:small}]
Let $u_{2}$ be an eigenfunction of the Laplacian operator corresponding to $\lambda_2(G,0)$ and assume that \eqref{eq:small} holds for a given leaf $v$.

We denote by $[v,v+\varepsilon]$ a section of length $\varepsilon>0$ of the incident edge $e$ adjacent to $v$, where $\varepsilon$ is chosen, first, such that $l_e > \varepsilon$, and second, such that
\begin{displaymath}
    u_{2}^{2}(y)<\frac{1}{L}\int_{G}u_{2}^{2}(x)\,\dx = \frac{1}{L}
\end{displaymath}
for all $y \in [v,v+\varepsilon]$, which is possible by \eqref{eq:small} and continuity of $u$. (Here we are also assuming $u_2$ is normalized to have $L^2$-norm $1$.)

We now define a piecewise linear, convex potential on $G$ with support on $[v,v+\varepsilon]$, by
\begin{displaymath}
    q(x):= \begin{cases} \varepsilon - \dist(x,v)
    \qquad &\text{on } [v,v+\varepsilon],\\ 0 \qquad 
    &\text{elsewhere};\end{cases}
\end{displaymath}
then $q$ is integrable, $q_t(x):=tq(x)$ is convex for all $t \geq 0$ and clearly differentiable in $t$ for all $t \in \R$, and $\|q_t\|_\infty = t\varepsilon < M$ for $t$ small enough; thus in particular $q_t \in C_{G,M}$ for $t \geq 0$ small.

Since the normalized eigenfunction $u_1(x) = \frac{1}{\sqrt{L}}$ is constant, it follows that
\begin{displaymath}
\begin{aligned}
    \int_G q(x) [u_2^2(x)-u_1^2(x)]\,\dx
    &=\int_{[v,v+\varepsilon]} q(x) [u_2^2(x)-u_1^2(x)]\,\dx\\
    &<\int_{[v,v+\varepsilon]} q(x)\left(\frac{1}{L}-\frac{1}{L}\right)\,\dx = 0.
\end{aligned}
\end{displaymath}
By Corollary~\ref{cor:feynman-hellman}, the zero potential does not minimize $\Gamma [q]$ in $C_{G,M}$.
\end{proof}

\begin{cor}
\label{cor:vanishing}
Given a compact tree $G$, if there is an eigenfunction of the Laplacian operator corresponding to $\lambda_2(G,0)$ which vanishes identically on an edge, then, for any $M>0$ the potential $q=0$ does not minimize $\Gamma [q]$ in $C_{G,M}$.
\end{cor}

\begin{proof}
Let $u_2$ be any such eigenfunction which vanishes identically on some edge. Since $G$ is a tree and the regions $G_+ := \overline{\{x \in G: u_2 (x)> 0\}}$ and $G_- := \overline{\{x \in G: u_2 (x)< 0\}}$ are necessarily connected (and have exactly one point of intersection $v_0$, without loss of generality a vertex), the zero set $G_0 := \{x \in G: u_2(x)=0\}$ must also be connected and intersect $G_+$ and $G_-$ at $v_0$ only.

In particular, $G_0$ is a subtree of $G$ of positive measure attached to $G\setminus G_0$ at $v_0$, and hence contains at least one leaf $v \in \mathcal{V}(G)$ of $G$. Since $u_2(v)=0$, Lemma~\ref{lem:small} yields the conclusion.
\end{proof}

\begin{cor}
\label{cor:multiple}
Given a compact tree $G$, if the multiplicity of $\lambda_2(G,0)$ is at least two, then for any $M>0$ the potential $q=0$ does not minimize $\Gamma [q]$ in $C_{G,M}$.
\end{cor}

\begin{proof}
Let $u_{2,1}$, $u_{2,2}$ be any two linearly independent eigenfunctions on $G$ associated with $\lambda_2(G,0)$ and let $v \in \mathcal{V}(G)$ be any leaf. Then there exist $a_1,a_2 \in \R$, not both zero, such that $a_1u_{2,1}(v)+a_2u_{2,2}(v) = 0$. Since $a_1u_{2,1}+a_2u_{2,2}$ is still a (nonzero) eigenfunction, the conclusion follows immediately from Lemma~\ref{lem:small}.
\end{proof}

Obviously, both corollaries allow the construction of various concrete examples for which the constant potential is not minimal; for example, any equilateral star graph on at least three edges satisfies the conditions of both. However, they also show that, given any compact tree $G$, it is possible to append an arbitrarily short edge to $G$ in such a way that on the resulting tree the constant potential is not minimal. Such a principle was also exploited in \cite{KeRo} (e.g., Theorem 4.5 there; see also \cite{James}), and reinforces the heuristic idea that in some analytic sense the set of degree-one vertices of a graph is not a ``good'' notion of boundary of the graph.

\begin{thm}\label{thm:constant-perturbation}
Let $G$ be a compact tree and suppose $\lambda_2(G,0)$ is simple. For any $\varepsilon>0$ there exists a tree $G_\varepsilon$ formed by attaching an edge of length $\varepsilon$ to $G$, such that for any $M>0$ the potential $q=0$ does not minimize $\Gamma [q]$ in $C_{G_\varepsilon,M}$.
\end{thm}

The theorem is an immediate consequence of the following lemma and Corollary~\ref{cor:vanishing}. While the conclusion should certainly still hold if $\lambda_2(G,0)$ is not simple, this case is not really of interest in light of Corollary~\ref{cor:multiple}.

\begin{lem}
\label{lem:graph-extension}
Let $G$ be a compact tree and suppose for some $q \in L^\infty(G)$ with $\|q\|_\infty \leq M$ that $\lambda_2 (G,q)$ is simple. Denote by $v_0$ the unique zero of the (any) corresponding eigenfunction $u_2$ and assume without loss of generality that $v_0 \in \mathcal{V}(G)$.

Then there exists $\varepsilon_{0}>0$ such that, for any $\varepsilon\in (0,\varepsilon_{0})$, if we attach a pendant edge $e_{\varepsilon}$ of length $\varepsilon$ to $G$ at $v_{0}$ to form a new graph $G_\varepsilon$, and extend $q$ to a potential $q_{\varepsilon}$ on $G_{\varepsilon}$ in such a way that $\|q_\varepsilon\|_\infty \leq M$, then $\lambda_{2}(G_{\varepsilon},q_{\varepsilon})$ is simple, and
$$\lambda_{2}(G,q)=\lambda_{2}(G_{\varepsilon},q_{\varepsilon}).
$$
Moreover,
 \begin{equation}\label{ep4}
u_{2,\varepsilon}(x)=\begin{cases} u_2(x) \qquad &\text{if } x \in G,\\
    0 \qquad &\text{if } x \in e_\varepsilon,\end{cases}
\end{equation}
is a corresponding eigenfunction.
\end{lem}

(We recall, again, that we are assuming standard vertex conditions throughout.)

\begin{proof}

\noindent For $\varepsilon>0$ we define $u_{2,\varepsilon}$ on $G_{\varepsilon}$ by \eqref{ep4}. Then $u_{2,\varepsilon}$ solves $-u''_{2,\varepsilon}+qu_{2,\varepsilon}=\lambda_{2}(G,q)u_{2,\varepsilon}$ pointwise on all edges of $G$, while it trivially satisfies $-u''_{2,\varepsilon}+q_\varepsilon u_{2,\varepsilon}=\lambda_{2}(G,q)u_{2,\varepsilon}$ on $e_\varepsilon$. Since it also satisfies standard conditions at all vertices, we conclude it is an eigenfunction on $G_\varepsilon$, and thus $\lambda_{2}(G,q)$ is \emph{some} eigenvalue of $G_\varepsilon$.

Now it is known (see \cite{Berkolaiko3}) that, for each fixed $k \geq 1$, there is eigenvalue convergence $\lambda_{k}(G_{\varepsilon},q_{\varepsilon})\rightarrow \lambda_{k}(G,q)$ as $\varepsilon \to 0$. Since by construction
\begin{displaymath}
    \lambda_{1}(G,q)<\lambda_{2}(G,q)<\lambda_{3}(G,q)
\end{displaymath}
and $\lambda_2 (G_\varepsilon, q_\varepsilon) = \lambda_{2}(G,q) \to \lambda_2 (G,q)$, eigenvalue convergence implies that for $\varepsilon>0$ small enough,
\begin{displaymath}
    \lambda_{1}(G_{\varepsilon},q_{\varepsilon})<\lambda_{2}(G_{\varepsilon},q_{\varepsilon}) = \lambda_2 (G,q) <\lambda_{3}(G_{\varepsilon},q_{\varepsilon}),
\end{displaymath}
and $u_{2,\varepsilon}$ is a second eigenfunction for $G_\varepsilon$.
\end{proof}

\subsection{Graphs with a long edge}
We now turn to the second kind of perturbation argument to which we alluded above: here we are interested in ``global'' perturbations of $q=0$ by a linear potential.

\begin{defi}\label{def:signed-distance}
Given a compact, connected graph $G$, let $x_{0}\in G$ be arbitrary, such that $G\setminus\{x_0\}$ is disjoint. Arbitrarily assign the respective labels $G^+$ and $G^-$ to two disjoint (and not necessarily connected) subsets of $G\setminus\{x_{0}\}$ whose union forms $G\setminus\{x_{0}\}$. We define the signed distance function $\sigma_{x_{0}}: G \to \R$ (with respect to the given labeling) by:
\begin{equation}
\label{ep3}
\sigma_{x_{0}}(x)=
\begin{cases}
    \dist (x,x_0) \qquad &\text{if } x \in G^+,\\
    -\dist (x,x_0) \qquad &\text{if } x \in G^-,\\
    0 \qquad &\text{if } x = x_0.
\end{cases}
\end{equation}
\end{defi}

Clearly $\sigma_{x_0}$ is always \emph{affine}, for any $x_0 \in G$; however it is not generally convex, cf.\ Remark~\ref{rem:affine}:

\begin{lem}
\label{lem:convex-sigma}
Let $G$ be a compact tree and $x_{0}\in G$. Then $\sigma_{x_{0}}(x)$ is convex if and only if $G^-$ is an interval.
\end{lem}

If $q$ is convex and $\sigma_{x_0}$ satisfies the conditions of the above lemma, then clearly $q+t\sigma_{x_0}$ is convex for all $t \geq 0$ (but not for $t<0$!).

We will give a single (numerical) example illustrating how $\sigma_{x_0}$ can be used to show that on a graph with a long edge $q=0$ may not minimize the spectral gap in any class of the form $C_{G,M}$.

\begin{ex}
Let $G$ be a star graph with $4$ vertices $\mathcal{V}=\left\lbrace v_{0},v_{1},v_{2},v_{4} \right\rbrace $ and $3$ edges $\mathcal{E}=\left\lbrace e_{1},e_{2},e_{4} \right\rbrace$ where $e_{i}=[v_{0},v_{i}]$ has length $i$ for all $i=1,2,4$; in particular, $v_0$ will be the central vertex of degree $3$.
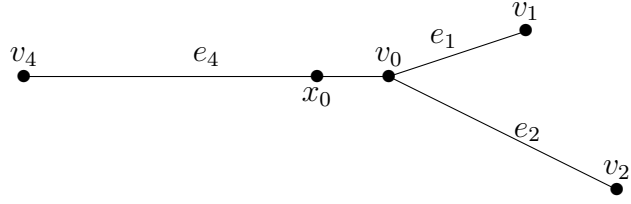
\begin{center}
\begin{tikzpicture}[scale=0.6]
\draw (-8,0) -- (0,0);
\draw (-4,0) node[above]{$e_{4}$} ;
\draw (1.2,0.4) node[above]{$e_{1}$} ;
\draw (2.5,-1.23) node[right]{$e_{2}$} ;
\draw (-8,0) node[above ]  { $v_{4}$}node{$\bullet$};
\draw (3,1) -- (0,0);
\draw (3,1) node[above ] {$v_{1}$} node{$\bullet$}; 
\draw[black] (0,0) -- (0,0) node[pos=0.5,above] {$v_{0}$};
\draw (5,-2.5) -- (0,0);
\draw (5,-2.5) node[above ] {$v_{2}$} node{$\bullet$}; 
\draw (0,0) node[below right] {} node{$\bullet$};
\draw (-1.57,0) node[below] {$x_0$} node{$\bullet$};
\end{tikzpicture}
\end{center}

We claim that in this case the constant potential $q=0$ is not optimal in $C_{G,M}$ for any $M>0$, something we will show numerically using Corollary~\ref{cor:feynman-hellman} and a potential of the form $t\sigma_{x_0}$, $t \geq 0$ as a perturbation; we observe that $t\sigma_{x_0} \in C_{G,M}$ for all $t \geq 0$ small enough (how small depending on $x_0$). We first find the point $x_0$ in the interior of $e_4$ such that
\begin{displaymath}
    \int_G \sigma_{x_0}(x)\,\dx = 0
\end{displaymath}
(and $G^-$ corresponds to the pendant part $[x_0,v_4]$ of $e_4$, so that Lemma~\ref{lem:convex-sigma} holds; the existence of such a $x_0$ is only guaranteed since $e_4$ is sufficiently long). Writing $e_i \sim [0,i]$, where $[0,i] \sim [v_0,v_i]$ for $i=1,2,4$, we wish to have
\begin{displaymath}
    \int_{G}\sigma_{x_{0}}(x)\,\dx=-\int_{0}^{4-x_{0}}x\,\dx
    +\int_{0}^{x_{0}}x\,\dx+\int_{x_{0}}^{x_{0}+1}2x\,\dx
    +\int_{x_{0}+1}^{x_{0}+2}x\,\dx=0.
\end{displaymath}
This leads to $x_0 = \frac{11}{14}$, that is, $x_0$ is at distance $\frac{11}{14}$ from $v_0$.

We now use $t\sigma_{x_0}$, $t \geq 0$, as a perturbation in Corollary~\ref{cor:feynman-hellman}. Since any first eigenfunction associated with $q=0$ is constant, condition \eqref{eq:feynman-hellman-condition} reduces to
\begin{equation}
\label{eq:3-star-long-edge-example}
    \int_G \sigma_{x_0}(x) u_2(x)^2\dx < 0
\end{equation}
for some second eigenfunction $u_2$ associated with $q=0$.

To show that \eqref{eq:3-star-long-edge-example} holds, we will take a different labeling convention more convenient for studying the eigenfunctions: we associate $e_i$, $i=1,2,3$, with $[0,i]$, where now $0 \sim v_i$ and $i \sim v_0$. This means that, on each edge, by the Neumann condition at $0$, any eigenfunction of the operator $H$ for $q=0$ on $e_{i}$ must be of the form $u(x)=A_{i}\cos(kx)$, $x \in [0,i]$, where $k=\sqrt{\lambda}$.

Now continuity and the Kirchhoff condition at $v_0$ lead to the transcendental equation
\begin{displaymath}
    \tan(k)+\tan(2k)+\tan(4k)=0
\end{displaymath}
(note that $\lambda_2$ will correspond to the smallest positive solution $k>0$) and hence to the choice
\begin{displaymath}
    A_1 = \frac{\cos (4k)}{\cos (k)}, \qquad
    A_2 = \frac{\cos (4k)}{\cos (2k)}, \qquad A_4 = 1,
\end{displaymath}
unique up to scalar multiples.

Using these values and the value of $x_0$ found above, with our second labeling convention the integral in \eqref{eq:3-star-long-edge-example} becomes
\begin{align*}
\int_{G}\sigma_{x_{0}}(x)u_2(x)^2\,\dx
&=\left[ \frac{\cos(4k)}{\cos(k)}\right]^{2}\int_0^1\left(\frac{25}{14}-x \right) \cos^{2}(kx)\,\dx\\
&+\left[ \frac{\cos(4k)}{\cos(2k)}\right]^{2} \int_0^2\left(\frac{39}{14}-x \right) \cos^{2}(kx)\,\dx\\
&+\displaystyle\int_0^4\left(x-\frac{45}{14} \right) \cos^{2}(kx) \,\dx.
\end{align*}
A numerical computation (performed using Mathematica) yields $k\approx 0.502642$ and thus
\begin{displaymath}
    \int_{G}\sigma_{x_{0}}(x)u_2(x)^2\,\dx\approx -1.46034,
\end{displaymath}
which is clearly negative (well outside the range of machine error). This establishes \eqref{eq:3-star-long-edge-example} numerically and hence shows that $q=0$ does not minimize the spectral gap in $C_{G,M}$.
\end{ex}

The above example suggests that, more generally, if $G$ consists of a relatively long edge to which a relatively small ``decoration'' is attached, then we do not expect the constant potential to be optimal: indeed, we can find a suitable $\sigma_{x_0}$ as a perturbation, and thus it is enough that \eqref{eq:3-star-long-edge-example} hold for $q=0$ not to be optimal.

\section{Bounds on the fundamental gap}
\label{sec:bounds}

In this section we will consider bounds on the fundamental gap $\Gamma[G,q]$ not just as a function of the potential $q$ but in particular as a function of the graph $G$, mostly in terms of the diameter $D$ of $G$.

Diameter is a well-known quantity on convex domains; along with convex potentials it is the ``other half'' of the fundamental gap conjecture proved by Andrews and Clutterbuck \cite{Andrews}: on a convex domain $\Omega \subset \R^n$ and with a suitable convex potential $q: \Omega \to \R$, the fundamental gap of the Schr\"odinger operator $-\Delta + q$ with Dirichlet conditions satisfies
\begin{displaymath}
    \Gamma[\Omega,q] \geq \frac{3\pi^2}{D^2}.
\end{displaymath}
For Neumann conditions (and $q=0$) the corresponding bound is due to Payne--Weinberger \cite{Payne}, $\Gamma[\Omega, 0] \geq \frac{\pi^2}{D^2}$. For the corresponding free Laplacian on metric graphs (i.e.\ $q=0$ and standard conditions) the situation is a bit more complicated: it was shown in \cite[Section~5]{KKMM} that, on a general compact graph, no upper or lower bound on $\Gamma[G,0]$ is possible only in terms of $D$; but if $G$ is a \emph{tree}, then Rohleder's theorem \cite[Theorem~3.4]{Rohleder} yields in particular that
\begin{displaymath}
    \Gamma[G,0] \leq \frac{\pi^2}{D^2},
\end{displaymath}
the opposite bound to the one of the Payne--Weinberger theorems.

Here we will give purely negative results in the case of lower bounds (still with standard vertex conditions throughout):
\begin{enumerate}
\item Even if $q=0$ and $G$ is a tree, no lower bound on $\Gamma[G,0]$ is possible in terms of diameter alone. This is a trivial modification of the example in \cite[Section~5.1]{KKMM}, whose details we omit (just replace the ``flower dumbbells'' there with ``star dumbbells'').
\item It is possible to find a \emph{fixed star graph} $G$ and a sequence of convex potentials $q_n$ (with $\|q_n\|_\infty \to \infty$), such that $\Gamma[G,q_n] \to 0$. See Example~\ref{ex:star-graph-no-bound-1}.
\item Alternatively, we can find a sequence of star graphs $G_n$ of fixed diameter $D$ and step potentials $q_n \in SW_{G_n,M}$ (for some constant $M > \frac{\pi^2}{D^2}$) such that $\Gamma [G_n,q_n] \to 0$. See Example~\ref{ex:star-graph-no-bound-2}.
\end{enumerate}

\begin{rem}
In terms of lower bounds, the only positive result is that on a \emph{fixed} compact tree $G$, in any of our given classes $C_{G,M}$, $SW_{G,M}$ etc., since there \emph{is} a minimizing potential $q^\ast$ we obtain the bound
\begin{displaymath}
    \Gamma[G,q] \geq \Gamma[G,q^\ast] > 0
\end{displaymath}
for all $q$ in the corresponding class; the point is that this bound is positive (since no potential can generate a first eigenvalue of multiplicity two). Thus on a fixed tree $G$ and within a fixed class we cannot find $q_n$ such that $\Gamma[G,q_n] \to 0$. But in the present work we will not study further the question of finding explicit estimates for given graphs and given classes.
\end{rem}

For upper bounds Rohleder's theorem can be trivially generalized to provide the bound
\begin{equation}
\label{eq:rohleder-with-q}
    \Gamma[G,q] \leq \frac{\pi^2}{D^2} + \|q\|_\infty
\end{equation}
for any compact tree $G$ and any potential $q \geq 0$: just use the Rayleigh quotient and \cite[Theorem~3.4]{Rohleder}, respectively, to obtain
\begin{displaymath}
    \lambda_2(G,q) \leq \lambda_2(G,\|q\|_\infty)
    = \lambda_2 (G,0) + \|q\|_\infty \leq \frac{\pi^2}{D^2}
    +\|q\|_\infty,
\end{displaymath}
while $\lambda_1(G,q) \geq 0$. However:
\begin{enumerate}
\item[(4)] There exist a sequence of trees $G_n$ of fixed diameter and single-well potentials $q_n$ on $G_n$ (with $\|q_n\|_\infty \to \infty$) for which $\Gamma[G_n,q_n] \to \infty$; the construction is rather complicated and we will only sketch it in Example~\ref{ex:star-graph-no-bound-3}. Notably, in this case the \emph{ratio} of the first two eigenvalues will still remain bounded.
\end{enumerate}

\begin{ex}
\label{ex:star-graph-no-bound-1}
Let $G$ be a $3$-star graph of three edges $e_-,e_+,e_\varepsilon$ connected at a central vertex $v_0$. We denote the other vertex incident with $e_\pm$ by $v_\pm$, and the other vertex incident with $e_\varepsilon$ by $v_\varepsilon$, respectively. We identify $e_\pm$ with the interval $[v_\pm,v_0]$, and $e_\varepsilon$ with the interval $[0,\varepsilon]$, where $0 \sim v_0$ and $\varepsilon \sim v_\varepsilon$.

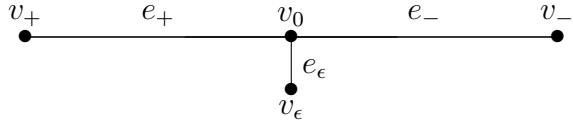
\begin{center}
\begin{tikzpicture}[scale=0.7]
\draw (-5,0) -- (5,0);
\draw (-2.5,0) node[above]{$e_{+}$} ;
\draw (2.5,0) node[above]{$e_{-}$} ;
\draw (0,-0.6) node[right]{$e_{\epsilon}$} ;

\draw (-5,0) node[above ]  { $v_{+}$}node{$\bullet$} -- (5,0) node[above ] {$v_{-}$}node{$\bullet$};
\draw[black] (-2,0) -- (2,0) node[pos=0.5,above] {$v_{0}$};

\draw (0,-1) -- (0,0);
\draw (0,-1) node[below ] {$v_{\epsilon}$} node{$\bullet$}; 

\draw (0,0) node[below right] {} node{$\bullet$};

\end{tikzpicture}
\end{center}

We normalize to take $|e_+|=|e_-|=1$, and we imagine $\varepsilon>0$ as being small but fixed; thus $D(G)=1$.

We consider the family of potentials
\begin{equation}
\label{eq:qt}
q_{t}(x)= \begin{cases} 0 \qquad &\text{if } x \in e_- \cup e_+,\\
    tx \qquad &\text{if } x \in e_\varepsilon \sim [0,\varepsilon],
    \end{cases}
\end{equation}
where $t \geq 0$. It is clear that $q_{t}$ is convex for all $t\geqslant0$, and $t\mapsto\lambda_{n}(t) := \lambda_n (G,q_t)$ is an increasing function of $t$ for each fixed $n \geq 1$.

We claim that, for $\varepsilon>0$ small but fixed, $\lambda_2(t) = \pi^2$ for all $t \geq 0$, while $\lambda_1(t) \to \pi^2$ (and thus $\Gamma[G,q_t] \to 0$) as $t \to \infty$.

We first observe that, for $\varepsilon < 1$ fixed, $\lambda_2(t) = \pi^2 = \lambda_2(e_-\cup e_+, 0)$ for all $t \geq 0$; and this eigenvalue is simple with eigenfunction supported entirely on $e_-\cup e_+$. Indeed, if $\varepsilon < 1$, then we see directly that $\lambda_2(0) = \pi^2$ (with eigenfunction identically zero on $e_\varepsilon$); we also see that this eigenfunction remains an eigenfunction for any $t>0$, and thus $\pi^2$ is always an eigenvalue. Since $\lambda_2(t)$ is a monotonically increasing function of $t$ and $\lambda_1(t)$ is always simple, the only possibility is that $\lambda_2(t) = \pi^2$ for all $t$. Simplicity of $\lambda_2(0)$ and monotonicity of $\lambda_3(t)$ imply that $\lambda_2(t)$ is always simple.

We next study $\lambda_1(t)$ and its associated positive eigenfunction $u_1(x):=u_1(t,x)$; by \cite[Theorem~1]{Kurasov}, $u_1(x)>0$ for all $x \in G$; moreover, by simplicity of $\lambda_1$, $u_1$ is symmetric on $e_-\cup e_+$.

\emph{Step 1.} We claim that, for sufficiently large fixed $t$, $u_1$ is monotonically decreasing on $e_\varepsilon$ from the central vertex $v_0 \sim 0$ to $v_\varepsilon \sim \varepsilon$.

First, note that $u_1'(\varepsilon)=0$ for all $t>0$ by the Kirchhoff condition; while for $x$ sufficiently close to $\varepsilon \sim v_\varepsilon$, $u_1$ is decreasing from $x$ towards $\varepsilon$, since
\begin{displaymath}
    -u_1'(x) = u_1'(\varepsilon) - u_1'(x) = \int_x^\varepsilon ts - \lambda_1(t) u_1(s)\,{\rm d}s > 0
\end{displaymath}
for $t$ large enough that $tx > \pi^2 \geq \lambda_1(t)$ (where we have used the eigenvalue equation to obtain the integral representation). Now if $u_1$ is not monotonic, then there exists $x_0 \in (0,\varepsilon)$ such that $u_1'(x_0)=0$, and hence $u_1$ is the first eigenfunction, and $\lambda_1(t)$ the first eigenvalue, of the interval $[x_0,\varepsilon]$ with Neumann conditions at the endpoints and the linear potential $q_t(x)=tx$. By a direct monotonicity argument, this is no smaller than the first eigenvalue of the interval $[0,\varepsilon]$ with potential $q_t(x)=tx$. But it is known that this latter eigenvalue diverges as $t \to \infty$ (see, e.g., \cite[Theorem~4.1]{ahrami}), a contradiction to $\lambda_1(t) \leq \pi^2$ for all $t \geq 0$.

\emph{Step 2.} Normalize $u_1$ to have $L^2$-norm $1$ on $G$; then $u_1(v_\varepsilon) \to 0$ as $t \to \infty$. To see this we first observe that, from the Rayleigh quotient and Step 1,
\begin{displaymath}
    \pi^2 \geq \lambda_1(t) > \int_0^\varepsilon q_t(x)u_1(x)^2\,\dx
    \geq u_1(v_\varepsilon)^2 \int_{0}^\varepsilon q_t(x)\,\dx.
\end{displaymath}
Since the latter integral diverges, we certainly have that $u_1(v_\varepsilon)^2 \to 0$.

\emph{Step 3.} Under the same normalization, $u_1(v_0) = \|u_1\|_{\infty,v_\varepsilon} \to 0$ as well. First observe that the argument of Step 2 implies that $u_1(x) \to 0$ for every \emph{fixed} $x \in e_\varepsilon$. If the claim of this step is not true, then there exist $\delta>0$ and a sequence $t_n \to \infty$ such that $u_1(v_0) \geq \delta$ for all $n \in \N$. Fix $x \in (0,\varepsilon) \subset v_\varepsilon$, to be chosen precisely later. Then, by the Fundamental Theorem of Calculus and Cauchy--Schwarz, for this $x$ and for \emph{any} $n \in \N$,
\begin{displaymath}
    u_1(v_0) = u_1(x) + \int_x^{v_0} u_1'(s)\,\textrm{d}s
    \leq u_1(x) + \|u_1'\|_{2,G}\sqrt{\dist (x,v_0)} < u_1(x) + \pi \sqrt{\dist (x,v_0)}.
\end{displaymath}
Choose $x$ close enough to $v_0$ such that $\pi \sqrt{\dist (x,v_0)} \leq \frac{\delta}{3}$, and then $n$ large enough that $u_1(x) \leq \frac{\delta}{3}$; then $u_1(v_0) \leq \frac{2\delta}{3}$ for all $n$ large enough, contradiction.

\emph{Step 4.} $\lambda_1 (t) \to \pi^2$ as $t \to \infty$. We observe that $u_1|_{e_- \cup e_+}$ is an eigenfunction, and $\lambda_2(t)$ an eigenvalue (necessarily the first, $\lambda_1^\alpha(e_-\cup e_+,0)$, since $u_1$ is positive) of the Laplacian on $e_- \cup e_+$ with standard conditions at $v_\pm$ and a Robin (delta) condition at $v_0$ of strength
\begin{displaymath}
    \alpha(t) = \frac{1}{u_1(v_0)}\sum_{x_i \in v} \partial u_1(x_i),
\end{displaymath}
where the sum is over the two edges $e_\pm$. (This is an example of ``cutting a graph along an eigenfunction'', as discussed in \cite[Section~3.1]{Berkolaiko-1}.) Since under our normalization $\|u_1\|_{2,G} = 1$ we have that $u_1(v_0) \to 0$, but clearly the derivative of $u_1$ on $e_- \cup e_+$ remains bounded away from zero, we conclude that $\alpha(t) \to \infty$ as $t \to \infty$, whence
\begin{displaymath}
    \lambda_1 (t) = \lambda_1^{\alpha(t)}(e_- \cup e_+,0)
    \to \lambda_1^\infty (e_- \cup e_+,0) = \pi^2,
\end{displaymath}
the first eigenvalue of the Laplacian on $e_- \cup e_+$ with standard (Neumann) conditions at $v_\pm$ and a Dirichlet condition at $v_0$.

This completes the proof that $\lambda_1 (G,q_t) \to \pi^2$, and hence $\Gamma [G,q_t] \to 0$.
\end{ex}

\begin{ex}
\label{ex:star-graph-no-bound-2}
We construct a variant of the above example: given $M > \pi^2$, a sequence of star graphs $G_{n}$ with $D(G_{n})=1$ for all $n\in \N$ and potentials $q_n \in SW_{G_n,M}$ (thus, in particular, $\|q_n\|_\infty$ remains bounded) such that
\begin{equation}
\label{eq:star-graph-no-bound-2}
    \lambda_{2}(G_{n},q_{n})=\pi^{2},\qquad \lambda_{1}(G_{n},q_{n})\to \pi^{2}.
\end{equation}
Fix $\varepsilon>0$ small and form $G_{n}$ by attaching two edges $e_{\pm}$ of length $\frac{1}{2}$ each and $n$ edges of length $\varepsilon$ at a central vertex $v_{0}$ of degree $n+2$. As in the previous example, we identify each shorter edge with $[0,\varepsilon]$, where $0$ corresponds to $v_0$, and each longer edge with $[0,\frac{1}{2}]$, where $0$ corresponds to $v_0$.

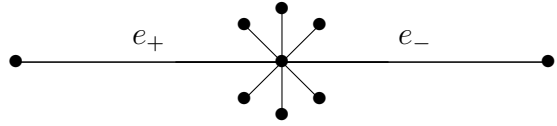
\begin{center}
\begin{tikzpicture}[scale=0.7]
\draw (-5,0) -- (5,0);
\draw (-2.5,0) node[above]{$e_{+}$} ;
\draw (2.5,0) node[above]{$e_{-}$} ;

\draw (-5,0) node[ below right ]  { }node{$\bullet$} -- (5,0) node[below right ] {}node{$\bullet$};
\draw[black] (-2,0) -- (2,0) node[pos=0.5,above] {};

\draw (0,-1) -- (0,0);
\draw (0,-1) node[below ] {} node{}; 

\draw (0,0) node[below right] {} node{$\bullet$};
\draw (0,-1) node[below right] {} node{$\bullet$};

\draw[black] (-2,0) -- (2,0) node[pos=0.5,above right] {};

\draw (0,1) -- (0,0);
\draw (0,1) node[below ] {} node{}; 

\draw (0,1) node[below right] {} node{$\bullet$};
node[pos=0.5,above] {};

\draw (0.707,-0.707) -- (0,0);
\draw (0.707,-0.707) node[below ] {} node{}; 

\draw (0.707,-0.707) node[below right] {} node{$\bullet$};
node[pos=0.5,above] {};

\draw (-0.707,0.707) -- (0,0);
\draw (-0.707,0.707) node[below ] {} node{}; 

\draw (-0.707,0.707) node[below right] {} node{$\bullet$};

node[pos=0.5,above] {};

\draw (0.707,0.707) -- (0,0);
\draw (0.707,0.707) node[below ] {} node{}; 

\draw (0.707,0.707) node[below right] {} node{$\bullet$};
node[pos=0.5,above] {};

\draw (-0.707,-0.707) -- (0,0);
\draw (-0.707,-0.707) node[below ] {} node{}; 

\draw (-0.707,-0.707) node[below right] {} node{$\bullet$};
\end{tikzpicture}
\end{center}
For simplicity, we define
\begin{equation*}
    q_n(x) := \begin{cases} 0 \qquad &\text{if } x \in e_- \cup e_+,\\
    M \qquad &\text{otherwise,}\end{cases}
\end{equation*}
so that, in particular, $q_n \in SW_{G_n,M}$, and claim that \eqref{eq:star-graph-no-bound-2} holds for this choice of $G_n$ and $q_n$. (However, the same construction should work if $q_n$ is chosen linear on each short edge to form a convex potential; here it is the increasingly large number of short edges close to $v_0$ that replaces the increasingly large potential of \eqref{eq:qt}.)

We start by observing that under these assumptions, an easy variant of Lemma~\ref{lem:graph-extension} shows that, for $\varepsilon>0$ small but fixed (independent of $n$), $\lambda_2 (G_n,q_n) = \lambda_2 (e_- \cup e_+,0) = \pi^2$, and there is a corresponding eigenfunction vanishing on all $n$ short edges.

Now denote by $\psi_n$ the positive eigenfunction corresponding to $\lambda_1 (G_n, q_n)$ with $L^2$-norm $1$; simplicity of the eigenvalue implies that $\psi_n$ is symmetric about $v_0$ on $e_+ \cap e_-$, and likewise invariant under permutation of the shorter edges.

\emph{Step 1.} We claim that, under the assumption that $M \geq \pi^2$ and our notational convention for the edges, $\psi_n$ is monotonically increasing from $0 \sim v_0$ to $\frac{1}{2}$ on the longer edges, and monotonically decreasing from $0 \sim v_0$ to $\varepsilon$ on the shorter edges. To see this, on the short edge, using the Kirchhoff condition $\psi_n'(\varepsilon)=0$ and the eigenfunction equation, we have
\begin{displaymath}
    \psi_n'(x) = -\int_x^{\varepsilon}\psi_n''(t)\,\dt
    =\int_x^{\varepsilon} \underbrace{[\lambda_1(G_n,q_n)-q_n(t)]}_{\leq 0} \psi_n(t)\,\dt \leq 0
\end{displaymath}
for all $x \in (0,\varepsilon)$. Now the Kirchhoff condition at $v_0$ plus symmetry implies that $\psi_n'(0)>0$ on each of the long edges. Since $\psi_n|_{e_+}$, being positive, necessarily corresponds to the first eigenvalue of the Laplacian on each of $e_\pm$ with a suitable Robin condition at $0$ and a Neumann condition at $\frac{1}{2}$, it must be monotonically increasing on $e_\pm$ from $0$ to $\frac{1}{2}$. This proves the claim.

\emph{Step 2.} Under the normalization $\|\psi_n\|_2 = 1$, if $M>\pi^2$, then we claim that
\begin{equation}
\label{eq:star-qn-step2}
    \int_{e_- \cup e_+} |\psi_n(x)|^2\,\textrm{d}x \not\to 0
\end{equation}
as $n \to \infty$. 
Denoting by $e_\varepsilon$ any of the $n$ shorter edges and using symmetry, we can estimate
\begin{displaymath}
\begin{aligned}
    \lambda_1 (G_n,q_n) &= \frac{n \int_{e_\varepsilon} |\psi_n'(x)|^2 + M |\psi_n(x)|^2\,\dx+\int_{e_- \cup e_+} |\psi_n'(x)|^2\,\dx}
    {n\int_{e_\varepsilon} |\psi_n(x)|^2\,\dx + \int_{e_ \cup e_+} |\psi_n(x)|^2\,\dx}\\
    & \geq \frac{n \int_{e_\varepsilon} |\psi_n'(x)|^2 + M |\psi_n(x)|^2\,\dx}
    {n\int_{e_\varepsilon} |\psi_n(x)|^2\,\dx + \int_{e_ \cup e_+} |\psi_n(x)|^2\,\dx}.
\end{aligned}
\end{displaymath}
In particular, if \eqref{eq:star-qn-step2} does \emph{not} hold, then
\begin{displaymath}
    \liminf_{n\to\infty}\lambda_1(G_n,q_n) \geq \liminf_{n\to\infty}
    \frac{\int_{e_\varepsilon} |\psi_n'(x)|^2 + M |\psi_n(x)|^2\,\dx}
    {\int_{e_\varepsilon} |\psi_n(x)|^2\,\dx} \geq M > \pi^2,
\end{displaymath}
a contradiction to $\lambda_1(G_n,q_n) < \lambda_2 (G_n,q_n) = \pi^2$.

\emph{Step 3.} We claim that, under the normalization $\|\psi_n\|_2=1$, $\psi_n(v_0) \to 0$. To see this, fix any short edge $e_\varepsilon \sim [0,\varepsilon]$; then by symmetry of $\psi_n$ and nonnegativity of $q_n$,
\begin{displaymath}
    \int_{e_\varepsilon} |\psi_n'(x)|^2\,\dx
    \leq \frac{1}{n}\int_{G_n} |\psi_n'(x)|^2\,\dx
    \leq \frac{1}{n}\lambda_1(G_n,q_n) \leq \frac{\pi^2}{n} \to 0.
\end{displaymath}
At the same time, using the monotonicity of $\psi_n$ established above,
\begin{displaymath}
    \psi_n (\varepsilon) = \min_{x \in e_\varepsilon} |\psi_n (x)| \to 0,
\end{displaymath}
as follows from $|G_n| \to \infty$ but $\|\psi_n\|_2 = 1$. Hence, by Cauchy--Schwarz,
\begin{displaymath}
    0 \leq \psi_n(0) - \underbrace{\psi_n(\varepsilon)}_{\to 0}
    = -\int_0^\varepsilon \psi_n'(x)\,\dx
    \leq \sqrt{\varepsilon} \left(\int_{e_\varepsilon}|\psi_n'(x)|^2\,\dx\right)^{1/2} \to 0,
\end{displaymath}
and so $\psi_n(0) = \psi_n(v_0) \to 0$ as well.

\emph{Step 4.} $\lambda_1 (G_n,q_n) \to \pi^2$ as $n \to \infty$. We first renormalize $\psi_n$ so that
\begin{displaymath}
    \int_{e_- \cup e_+} |\psi_n(x)|^2\,\dx = 1;
\end{displaymath}
by Steps 2 and 3, under this alternative normalization we still have that $\psi_n(v_0) \to 0$. Now as we have seen, restricted to $e_- \cup e_+$, $\psi_n$ is the first eigenfunction, and $\lambda_1 (G_n,q_n)$ the first eigenvalue, of the Laplacian on $e_- \cup e_+$ with a Robin condition, say of strength $\alpha_n > 0$, at $v_0$, and Neumann conditions at the other vertices. Now the fact that $\psi_n(v_0) \to 0$ implies that $\alpha_n \to \infty$, since otherwise
\begin{displaymath}
    \alpha_n = \frac{1}{\psi_n(v_0)}[\partial_\nu \psi_n|_{e_-}(v_0) + \partial_\nu \psi_n|_{e_+}(v_0)]
\end{displaymath}
would remain bounded as $n\to\infty$ (since the normal derivatives of $\psi_n$ as a Robin eigenfunction on an interval remain bounded under the chosen normalization). Continuity of $\lambda_1$ with respect to $\alpha_n$, including as $\alpha_n \to \infty$ (corresponding to a Dirichlet condition at $v_0$) implies that $\lambda_1 (G_n,q_n)$ converges to the first eigenvalue of the Laplacian on $e_- \cup e_+$ with a Dirichlet condition at $v_0$ and a Neumann condition at the other vertices; this eigenvalue is $\pi^2$. This completes the proof.
\end{ex}

\begin{ex}
\label{ex:star-graph-no-bound-3}
Here we will \emph{sketch} the construction of a family of tree graphs $G_n$ of fixed diameter $D=1$, with single-well potentials $q_n$, such that $\Gamma [G_n,q_n] \to \infty$. In this case both $\lambda_1 (G_n,q_n)$ and $\lambda_2 (G_n,q_n)$ will diverge to $\infty$. But at least in this example their \emph{ratio} will remain bounded; it converges to $4$ (the corresponding ratio for the Dirichlet Laplacian on an interval) as $n \to \infty$.

We take an interval of length $1$ and, for $n$ fixed, decorate it with $n$ very small stars equally spaced with distance $\frac{1}{n}$ between neighbors (each star should have the same number $m(n)$ of very short edges of equal length $\delta(n)<\!\!<\frac{1}{n}$, with $m$ and $\delta$ to be chosen more precisely later), and half that distance between the outermost stars and the terminal vertices of the original interval.
\begin{figure}[ht]
\centering
\begin{tikzpicture}
\draw[thick] (-5.25,0) -- (5.25,0);
\foreach \number in {-4.5,-3,-1.5,0,1.5,3,4.5}{
\draw[thick] (\number,-0.3) -- (\number,0.3);
\draw[thick] (\number-0.2,-0.2) -- (\number+0.2,0.2);
\draw[thick] (\number-0.2,0.2) -- (\number+0.2,-0.2);
\draw[fill] (\number,-0.3) circle (0.75pt);
\draw[fill] (\number,0.3) circle (0.75pt);
\draw[fill] (\number-0.2,-0.2) circle (0.75pt);
\draw[fill] (\number+0.2,-0.2) circle (0.75pt);
\draw[fill] (\number-0.2,0.2) circle (0.75pt);
\draw[fill] (\number+0.2,0.2) circle (0.75pt);
}
\draw[fill] (-5.25,0) circle (1.25pt);
\draw[fill] (5.25,0) circle (1.25pt);
\node at (-5.25,0) [anchor=north] {$-\tfrac{1}{2}$};
\node at (5.25,0) [anchor=north] {$\tfrac{1}{2}$};
\end{tikzpicture}
\label{fig:interval-with-decorations}
\end{figure}

For simplicity we take $q_n$ to be single well, of the form $q_n = M(n) >\!\!> 0$ on each star decoration, and $q_n = 0$ on the interval. (A convex potential similar to \eqref{eq:qt} should also work.)

Keeping $n$ fixed, and taking $m(n)$ and $M(n)$ sufficiently large and $\delta(n)$ sufficiently small, arguments similar to the ones used in Examples~\ref{ex:star-graph-no-bound-1} and~\ref{ex:star-graph-no-bound-2} imply that $\lambda_1 (G_n,q_n)$ and $\lambda_2 (G_n,q_n)$ are, respectively, arbitrarily close (within some pre-specified $\varepsilon>0$) of the first two eigenvalues of the Laplacian on an interval with Dirichlet conditions at each decoration point.
\begin{figure}[ht]
\centering
\begin{tikzpicture}
\draw[thick] (-5.25,0) -- (5.25,0);
\draw[fill] (-5.25,0) circle (1.25pt);
\draw[fill] (5.25,0) circle (1.25pt);
\foreach \number in {-4.5,-3,-1.5,0,1.5,3,4.5}{
\filldraw[color=white] (\number,0) circle (1.75pt);
\draw[thick] (\number,0) circle (1.75pt);
}
\end{tikzpicture}
\end{figure}

By construction, this interval has first eigenvalue $\frac{\pi^2}{n^2}$ and second eigenvalue $\frac{4\pi^2}{n^2}$; thus, for our given choice of $m(n)$, $M(n)$ and $\delta(n)$, $|\Gamma (G_n,q_n) - \frac{3\pi^2}{n^2}| < 2\varepsilon$. Thus $\Gamma (G_n,q_n) \to \infty$, although the ratio $\frac{\lambda_2(G_n,q_n)}{\lambda_1(G_n,q_n)}$ converges to $4$.
\end{ex}

\end{document}